\documentclass[10pt,a4paper]{article}
\usepackage[utf8]{inputenc}
\usepackage{amsmath}
\usepackage{amsfonts}
\usepackage{amssymb}
\usepackage{amsthm}
\usepackage{graphicx}
\usepackage{bussproofs}
\usepackage{multicol, xspace}

\newtheorem{theorem}{Theorem}
\newtheorem{lemma}{Lemma}
\newtheorem{corollary}{Corollary}

\theoremstyle{definition}
\newtheorem{definition}{Definition}
\newtheorem{remark}{Remark}

\newcommand{\isig}[1]{{\ensuremath {\mathrm{I}\Sigma_{#1}}}\xspace}

\newcommand{\X}{\mathcal{X}}
\newcommand{ \lab}{\mathrm{lab}}
\newcommand{\B}{\mathcal{B}}
\newcommand{\C}{\mathcal{C}}
\newcommand{\IL}{\mathbf{IL}}
\renewcommand{\L}{\mathbf{L}}
\newcommand{\GL}{\mathbf{GL}}
\newcommand{\R}{\mathbf{R}}
\renewcommand{\S}{\mathbf{S}}
\newcommand{\N}{\mathbb{N}}
\renewcommand{\P}{\mathsf{P}}

\newcommand{\M}{\mathsf{M}}
\newcommand{\F}{\mathcal{F}}
\newcommand{\subjj}{\mathrm{sub}}

\title{Labelled tableaux for interpretability logics}
\author{Tuomas Hakoniemi \and Joost J. Joosten}

\begin{document}
	\maketitle
	
	\paragraph{Abstract}
\noindent
In is paper we present a labelled tableau proof system that serves a
wide class of interpretability logics. The system is proved sound and
complete for any interpretability logic characterised by a frame
condition given by a set of universal strict first order Horn
sentences. As such, the current paper adds to a better proof-theoretical understanding of interpretability logics.

\section{Introduction}
	
	Provability logics like the G\"odel-L\"ob logic $\GL$ describe the structural behaviour of formalized provability in a simple propositional modal language. Interpretability logics are natural extensions of provability logics: they describe the structural behaviour of relative interpretability. 
	
	Essentially since Solovay's landmark paper \cite{Solovay:1976} we know that any $\Sigma_1$ sound theory that extends elementary arithmetic has the same provability logic $\GL$. The situation is very different for interpretability logics. Basically, for two different kind of theories we know the corresponding interpretability logics. 
	
	On the one hand Shavrukov \cite{Shavrukov:1988:InterpretabilityLogicPA} and independently Berarducci \cite{Berarducci:1990:InterpretabilityLogicPA} determined the interpretability logic of any sound and essentially reflexive theory like Peano Arithmetic to be $\IL\M$. 
 On the other hand, Visser has proven in \cite{Visser:1990:InterpretabilityLogic} that the interpretability logic of any sound and finitely axiomatised theory that proves the totality of super-exponentiation --like \isig{1}-- is $\IL\P$. 
	
	In case the base theory is neither finitely axiomatizable nor essentially reflexive, the situation turns out to be much more difficult and actually, to determine the interpretability logics in those situations remain open problems. Some partial results are known in the case of Primitive Recursive Arithmetic (\cite{BilkovaJonghJoosten:2009:PRA}) or in the case when we consider those modal principles that are provable in any reasonable arithmetical theory \cite{JoostenVisser:2000:IntLogicAll, GorisJoosten:2011:ANewPrinciple, Joosten:2015:TwoSeries}. In this sense, interpretability logics are in need of more study compared to provability logics.\\
\bigskip
	
	\emph{A miracle happens} writes Albert Visser as the first line of \cite{Visser:1997:OverviewIL}: Whereas provability is a $\Sigma_1$ complete predicate, the logic $\GL$ that governs its structural behaviour is nice, well behaved and simple\footnote{``Only" of complexity PSPACE.}. The situation with interpretability seems even more extreme since Shavrukov has shown in \cite{Shavrukov:1997:ReflexiveInfinitelyManyAxioms} that interpretability is $\Sigma_3$ complete and again, the modal logic describing the structural behaviour is nice. 
	
	But again here we see a discrepancy between provability and interpretability logics. In the case of interpretability logics we actually know to a much lesser extent \emph{how} nice they are. In particular, apart from some observations on the closed fragment (\cite{HajekSvejdar:1991:ClosedFormulasInterpretability, BouJoosten:2011}), close to nothing is known about the computational complexity of interpretability logics. Also, very little is known about well-behaved proof systems for interpretability logics with the sole exception of some work by Sasaki such as \cite{Sasaki:2002:CutFreeIL}. The current paper is intended to add to the proof-theoretic understanding of interpretability logics by studying labelled tableaux proof systems for them. 
	
	Tableaux proof systems are tightly related to sequent proof systems and are dual to them in many aspects. Rules in sequent proof systems typically have possibly multiple antecedents and single conclusions/succedents. Moreover, sequent-style proofs generally are trees that have the root at the bottom and are based on validity: the consequence of a rule is valid if (and often only if) all of the antecedents are valid. 
	
	On the other hand, rules in tableaux systems typically have single antecedent and possibly multiple succedents. Moreover, Tableaux proofs generally are trees that have the root at the top and are based on satisfiability: the antecedent of a rule is satisfiable if and only if some of the `consequents'/succedents is satisfiable.
	
	Labelled tableaux introduce extra devices to the syntax that aim to represent the accessibility relation in the corresponding Kripke-style semantics. This extra syntax allows us to give tableaux proof systems for many logics lacking a traditional one, where nodes of the tableaux carry only (sets of) formulas. A precursor for this idea of bringing a bit of semantics into the syntax appears already in \cite{Kanger:1957:ProvabilityInLogic}, and labelled tableaux as they are now known were introduced prominently by Fitting in \cite{Fitting:1972}. Standard references here are \cite{Fitting:1983:ProofMethods} and \cite{Gore:1999}. For more on the history and development of tableaux systems for modal logics see e.g. \cite{Gore:1999}.
	Naturally, labels have also been incorporated into sequent calculi. We refer the reader to \cite{Negri:2005:ProofAnalysisInModalLogic} for details on labelled sequent calculi for modal logics.
	\\
\bigskip
	
	{\bf Outline of the paper}. After introducing the necessary preliminaries in Section \ref{section:preliminaries}, we use Section \ref{section:TableauxForHornILs} to introduce the \emph{labelled tableaux system} for all interpretability logics $\IL{\sf X}$ characterised by a set of first order Horn formulas. It is shown how a \emph{systematic tableau} can be assigned to a finite set $\Gamma$ of formulas so that the tableau contains all the necessary information as to the satisfiability of $\Gamma$.
	
	In Sections \ref{section:soundness} and \ref{section:completeness} we show that the tableau proofs are sound and complete with respect to $\IL{\sf X}$-validity. In the last section we remark that our results concern most of the interpretability logics encountered in the literature, but not all.
	
\section{Preliminaries}\label{section:preliminaries}
	
	Interpretability logics are propositional modal logics with a unary modality $\Box$ whose dual modality $\Diamond$ is defined as $\Diamond := \neg \Box \neg$ corresponding to provability and consistency respectively, and a binary modality $\rhd$ corresponding to relative interpretability.
	
	In this paper we shall work with the Boolean connectives $\neg$ and $\to$. Thus, with $\sf Prop$ a countable set of propositional variables, the formulas $\mathcal F$ of interpretability logic are defined as 
	$$\mathcal F := {\sf Prop} \mid (\neg \mathcal F) \mid (\mathcal F \to \mathcal F ) \mid (\Box \mathcal F ) \mid (\mathcal F \rhd \mathcal F ).$$
	
	As always we will use the other connectives and Boolean constants freely since they can be defined from $\neg$ and $\to$. In order to use less parentheses we omit outer parentheses and shall say that $\neg, \Box$ and $\Diamond$ bind strongest, followed by the equally strong binding $\vee$ and $\wedge$ who bind stronger than $\rhd$ which in turn binds stronger than $\to$. Thus, for example, 
	\[
	p\rhd q \to p\wedge \Box r \rhd q\wedge \Box r
	\] 
	is short for 
	\[
	\Big[ (p\rhd q) \to \Big( \big( p\wedge (\Box r)\big) \rhd \big(q\wedge (\Box r) \big) \Big) \ \Big].
	\]
	
	\begin{definition}
		The axioms of the basic interpretability logic $\IL$ are, apart from all substitution instances (in the language of interpretability) of all propositional tautologies, given by the following axiom schemata
		\begin{itemize}
			\item[$L1$] 
			\ \  \ \ \ $\Box (A\to B) \to (\Box A \to \Box B)$;
			
			\item[$L2$] 
			\ \  \ \ \ $\Box (\Box A \to A) \to \Box A$;
			
			\item[$J1$]
			\ \  \ \ \ $\Box (A \to B) \to A\rhd B$;
			
			\item[$J2$]
			\ \  \ \ \ $(A \rhd B) \wedge (B\rhd C) \to A\rhd C$;
			
			\item[$J3$]
			\ \  \ \ \ $(A\rhd C) \wedge (B\rhd C) \to A\vee B \rhd C$;
			
			\item[$J4$]
			\ \  \ \ \ $A\rhd B \to (\Diamond A \to \Diamond B)$;
			
			\item[$J5$]
			\ \  \ \ \ $\Diamond A \rhd A$.
		\end{itemize}
		The rules are Modus Ponens and Necessitation: $A/ \Box A$.
	\end{definition}
	
	The following lemma collects two easily obtainable and well-known properties of $\IL$ that will play prominent role in our tableaux systems.
	\begin{lemma}\label{theorem:basicPropertiesIL}\ 
		\begin{enumerate}
			\item\label{theorem:basicPropertiesIL:Loeb}
			$\IL \vdash \ \Diamond A \to \Diamond (A \wedge \Box \neg A)$;
			
			\item\label{theorem:basicPropertiesIL:LoebILversion}
			$\IL \vdash  B\vartriangleright B \wedge \Box \neg B$.
		\end{enumerate}
	\end{lemma}

	The logic $\GL$ is the fragment of $\IL$ where the modal language is restricted to $\Box$. We shall consider various extensions of $\IL$. By $\IL{\sf X}$ we denote the logic that arises by adding the axiom scheme(s) $\sf X$ to $\IL$. The extensions of $\IL$ obtained by the following axiom schemes play a prominent role in the literature.
	
	\begin{itemize}
		
		
		\item[$\P$:]
		$A \rhd B \ \to \ \Box (A \rhd B)$;
		
		\item[$\M$:]
		$A \rhd B \ \to \ A \wedge \Box C \rhd B\wedge \Box  C$.
	\end{itemize}
	
	Interpretability logics allow for a relational semantics very much in the sense as $\GL$ does. 
	
	\begin{definition}
		An \emph{$\IL$-frame} is a triple $\langle W,R,S\rangle$ where $W$ is a non-empty domain set, whose members are often called \emph{worlds}, and $R$ is a binary relation on $W$ that is transitive and Noetherian (no infinite chains $x_0Rx_1Rx_2\ldots$). $S$ is a ternary relation on $W$ that is often considered as a collection $\{ S_x \}_{x\in W}$ of binary relations by fixing the first argument $x$ of the ternary $S$. It is required that each $S_x$ is a transitive and reflexive binary relation on $\{ y\in W \colon xRy\}$ satisfying the following property:
		$$\text{if } xRyRz\text{, then } yS_x z.$$
		
		An \emph{$\IL$-model} is a quadruple $\langle W,R,S,V \rangle$, where $\langle W,R,S\rangle$ is an $\IL$-frame and $V$ is a function assigning a collection $V(p)$ of worlds to a propositional variable $p$. Given an $\IL$-model $\langle W,R,S,V \rangle$ we define a forcing relation $\Vdash$ between worlds and formulas as usual:
		\begin{itemize}
			\item $M,x\Vdash p \Leftrightarrow x\in V(p)$;\\
			\item $M,x\Vdash \neg A \Leftrightarrow M,x \nVdash A$; \\
			\item $M,x\Vdash A \to B \Leftrightarrow M,x \nVdash A \text{ or } M,x\Vdash B$;  \\
			\item $M,x\Vdash \Box A \Leftrightarrow \forall y (xRy \Rightarrow M,y \Vdash A)$;\\
			\item $M,x \Vdash A\rhd B \Leftrightarrow \forall y \big(xRy \wedge \,M,y\Vdash A \Rightarrow \exists z (yS_xz \wedge \,M,z \Vdash B) \big)$.\\
		\end{itemize}
	\end{definition}
	
	We shall write $x\in M$ whenever $M = \langle W, R, S, V\rangle$ with $x\in W$ and likewise for frames. We write $M\vDash A$ to denote that $M, x\Vdash A$ for all $x\in M$. The above defined semantics is good in that one can prove completeness for $\IL$ as was first done in \cite{JonghVeltman:1990:ProvabilityLogicsForRelativeInterpretability}: \ 
	\[
	\IL \vdash A \ \Leftrightarrow \ \forall M\   M\vDash A. 
	\]
	
	An extension $\L$ of $\IL$ can be specified either axiomatically or semantically by restricting the class of models for example by specifying so called \emph{frame conditions}. We say that a frame $F := \langle W, R, S\rangle$ validates $A$ and we write $F\vDash A$ whenever for all valuations $V$ on $F$ we have $\langle W, R, S,V\rangle \vDash A$. 
	
	A set of first or higher order sentences $\mathcal{C}$ in the language $\{ \overline R, \overline S\}$ with $\overline R$ a binary and $\overline S$ a ternary first-order relation symbol is called a \emph{frame condition} for a logic $\L$ extending $\IL$ whenever we have 
	\[
	\langle W,R, S\rangle \vDash \L \ \Longleftrightarrow \ \langle W,R, S\rangle \vDash_{\sf fo/ho} \mathcal{C},
	\]
where in the right-hand side the interpretations of $\overline{R}$ and $\overline{S}$ are $R$ and $S$, respectively. Then we also say that $\mathcal{C}$ characterises the logic $\L$. As always, $\langle W,R, S\rangle \vDash \L$ denotes that $\langle W,R, S\rangle \vDash A$ for any theorem $A$ of $\L$ and we use a similar convention for models. From now on we will use the same symbol $R$ for $\overline R$ and its interpretation and likewise for $S$. 
	
	In case an axiomatic extension $\IL{\sf X}$ of $\IL$ is characterised by a set of strict universal Horn sentences in the language $\{R,S\}$ we say that $\IL{\sf X}$ is a \emph{Horn} logic. By a strict universal Horn sentence we mean a first order formula of the form
	\[
	\forall\ldots\forall\ (\varphi_1\wedge\ldots\wedge\varphi_n\rightarrow\psi),
	\]
	where $n\geq 0$, $\varphi_1,\ldots,\varphi_n$ and $\psi$ are atomic formulas and $\forall\ldots\forall$ denotes the universal closure. In case that $\IL{\sf X}$ is a Horn logic, we shall denote the corresponding frame condition by $\mathcal{C}_{\sf X}$ and call an $\IL$-frame satisfying $\mathcal{C}_{\sf X}$ an $\IL{\sf X}$-frame. 
	
	For example, the logic $\IL\P$ is characterised by the (universal closure of the) first order formula 
	\[
	xRy \wedge yRz \wedge zS_x u \ \to \ zS_yu,
	\]
	and $\IL$ is characterised by the empty frame condition (or $\forall x\ (xRx \to xRx)$ for that matter).
	
	In this paper, we shall -- given a frame -- reduce the binary modality $\rhd$ to a series of unary ones. We will do so, so that the corresponding tableaux rules become more amenable. Thus, given an $\IL$-frame $\F = \langle W,R,S\rangle$, we introduce new unary modal operators $\Box_x$ for each $x\in W$ and give the following truth definition for the operators in a model $M$ on the frame $\F$
	$$M,y\Vdash \Box_{x}A\Leftrightarrow \forall z\, ( yS_{x}z\Rightarrow M,z\Vdash A).$$
	
	Now it is easy to verify that for any $\IL$-model 
	\begin{equation}\label{translation}
	M,x\Vdash A\vartriangleright B \ \ \Leftrightarrow \ \ M,x\Vdash\Box(A\Rightarrow\neg\Box_x\neg B).
	\end{equation}
	
\section{Tableaux for Horn interpretability logics}\label{section:TableauxForHornILs}
	
	In this section we define a tableau proof method for interpretability logics which are Horn. Moreover, we will give a systematic tableau procedure for such $\IL{\sf X}$ that yields a canonical tableau given a finite set of formulas.
	
	As always, our tableaux will be downward growing trees. Each node of the tree carries a labelled formula. A labelled formula is a pair with a label and a formula. The label corresponds to a possible world where the formula is to be satisfied. We will show the unsatisfiability of a finite set of formulas in case all branches in the systematic tableau close (precise definition follow). In case the systematic tableau contains an open branch, that branch will carry information about a satisfying model.
	
	\begin{definition}
		Labels are strings composed of non-negative integers and letters R and S. The set of all labels is defined recursively as follows:
		\begin{itemize}
			\item $0$ is a label;
			\item If $\sigma$ is a label, then $\sigma R n$ is a label for all $n\in\N$;
			\item If $\sigma$ and $\rho$ are labels and $\rho$ is a strict non-empty prefix of $\sigma$, then $\sigma S_\rho n$ is a label for all $n\in\N$. 
		\end{itemize}
	\end{definition}

Now that we have a sufficiently large set of labels we will describe how we generically build (almost) frames from them.
	
	\begin{definition}[$\IL{\sf X}$-label structure]\label{definition:ILrelationsOnLabels}
		Given a Horn logic $\IL{\sf X}$ and a set of labels $\Lambda$, we define relations $\R_{\IL{\sf X}}^{\Lambda}$ and $\S_{\IL{\sf X}}^\Lambda$ on the set $\Lambda$ as the least relations on $\Lambda$ such that:\\
		\begin{enumerate}
			\item \label{definition:ILrelationsOnLabels:FirstItem}
			If $\sigma,\sigma R n\in\Lambda$, then $\langle \sigma, \sigma R n\rangle \in\R_{\IL{\sf X}}^\Lambda$ for all labels $\sigma$ and $n\in\N$;\\
			
			\item 
			If $\langle \sigma,\tau\rangle \in\R_{\IL{\sf X}}^\Lambda$ and $\langle \tau,\rho\rangle \in\R_{\IL{\sf X}}^\Lambda$, then $\langle \sigma,\rho\rangle \in\R_{\IL{\sf X}}^\Lambda$;\\
			
			\item 
			If $\sigma,\rho,\sigma S_\rho n\in\Lambda$, then $\langle \rho, \sigma, \sigma S_\rho n\rangle \in\S_{\IL{\sf X}}^\Lambda$ for all labels $\sigma$ and $\rho$ and all $n\in\N$;\\
			
			\item 
			If $\langle \sigma,\tau\rangle \in\R_{\IL{\sf X}}^\Lambda$, then $\langle \sigma,\tau,\tau\rangle \in\S_{\IL{\sf X}}^\Lambda$;\\
			
			\item 
			If $\langle \rho,\sigma\rangle \in\R_{\IL{\sf X}}^\Lambda$ and $\langle \sigma,\tau\rangle \in\R_{\IL{\sf X}}^\Lambda$, then $\langle \rho,\sigma,\tau\rangle \in\S_{\IL{\sf X}}^\Lambda$;\\
			
			\item 
			If $\langle \rho,\sigma,\tau\rangle \in\S_{\IL{\sf X}}^\Lambda$ and $\langle \rho,\tau,\upsilon\rangle \in\S_{\IL{\sf X}}^\Lambda$, then $\langle \rho,\sigma,\upsilon\rangle \in\S_{\IL{\sf X}}^\Lambda$;\\
			
			\item 
			If $\langle \rho,\sigma,\tau\rangle \in\S_{\IL{\sf X}}^\Lambda$, then $\langle \rho,\sigma\rangle \in\R_{\IL{\sf X}}^\Lambda$ and $\langle \rho,\tau\rangle \in\R_{\IL{\sf X}}^\Lambda$;\\
			
			\item \label{definition:ILrelationsOnLabels:LastItem}
			$\langle \Lambda,\R_{\IL{\sf X}}^\Lambda,\S_{\IL{\sf X}}^\Lambda\rangle \vDash_{\sf fo} \mathcal{C}_{\sf X}$.\\
			
		\end{enumerate}
	\end{definition}
	Note that the least relations exist since $\mathcal{C}_{\sf X}$ is a set of strict first order Horn sentences. If the context allows us so, we will drop both the sub- and the superscripts in $\R_{\IL{\sf X}}^\Lambda$ and $\S_{\IL{\sf X}}^\Lambda$. Moreover, when $\langle \rho,\sigma,\tau\rangle\in\S$ we will denote this by $\sigma\S_\rho\tau$ and likewise for $\R$. 
	
	Note that $\R$ is irreflexive in case $\IL{\sf X}$ is consistent and $\Lambda$ sufficiently nice. Moreover, apart from $\R$ being Noetherian, all the other properties of $\IL{\sf X}$-frames are satisfied: $\R$ is transitive; $\S_\rho$ is a relation on $\{ \sigma \in \Lambda \mid \rho \R \sigma \}$ that is transitive and reflexive so that $\rho \R \sigma \R \tau \Rightarrow \sigma \S_\rho \tau$. 
	
	We will now define the generating rules for tableaux for Horn interpretability logics. As mentioned, the nodes carry \emph{labelled formulas} which consist of a pair $\sigma :: A$, where $\sigma$ is a label and $A$ is a formula. Recall that the idea of the labels is, that they will correspond to worlds in a model where the corresponding formula will be satisfied if satisfiable. 
	
	The rules that we present are not entirely local since, for example, we have to guarantee that new labels have not yet been used in relevant parts of the tableau so far. Thus, we define the rules relative to a set of labels. 
	
	\begin{definition}[Tableau rules]
		Let $\IL{\sf X}$ be a Horn interpretability logic and let $\Lambda$ be a set of labels. The $\IL{\sf X}$-tableau rules with respect to $\Lambda$ are as follows: 
		
		\ \\ \
		\noindent Propositional rules:
		\begin{prooftree}
			\AxiomC{$\sigma :: \neg\neg A$}
			\LeftLabel{($\neg$)}\RightLabel{;}
			\UnaryInfC{$\sigma :: A$}
		\end{prooftree}
				\ \\ \
		\begin{prooftree}
			\AxiomC{$\sigma :: A\rightarrow B$}
			\LeftLabel{($\rightarrow$)}\RightLabel{;}
			\UnaryInfC{$\sigma :: \neg A \ \ \ \  \mid  \ \ \ \  \sigma :: B$}
		\end{prooftree}
				\ \\ \
		\begin{prooftree}
			\AxiomC{$\sigma :: \neg (A\rightarrow B)$}
			\LeftLabel{($\neg\rightarrow$)}\RightLabel{.}
			\UnaryInfC{$\sigma :: A$}
			\noLine
			\UnaryInfC{$\sigma :: \neg B$}
		\end{prooftree}
				\ \\ \
		\noindent $(\nu)$-rules:
		\begin{prooftree}
			\AxiomC{$\sigma :: \Box A$}
			\LeftLabel{$(\nu_\Box, \Lambda)$}
			\RightLabel{, when $\sigma\R\tau$;}
			\UnaryInfC{$\tau :: A$}
		\end{prooftree}
				\ \\ \
		\begin{prooftree}
			\AxiomC{$\sigma :: \Box_\rho A$}
			\LeftLabel{$(\nu_S, \Lambda)$}
			\RightLabel{, when $\sigma \S_{\rho}\tau$;}
			\UnaryInfC{$\tau :: A$}
		\end{prooftree}
				\ \\ \
		\begin{prooftree}
			\AxiomC{$\sigma :: A\vartriangleright B$}
			\LeftLabel{$(\nu_\vartriangleright, \Lambda)$}
			\RightLabel{, when $\sigma\R\tau$.}
			\UnaryInfC{$\tau :: \neg A\ \ \ \ \mid \ \ \ \  \tau :: \neg\Box_\sigma\neg B$}
		\end{prooftree}
				\ \\ \
		\noindent$(\pi)$-rules:
		\begin{prooftree}
			\AxiomC{$\sigma :: \neg\Box A$}
			\LeftLabel{$(\pi_\Box, \Lambda)$}
			\RightLabel{, where $n\in\N$ is such that $\sigma Rn\notin\Lambda$;}
			\UnaryInfC{$\sigma R n :: \neg A$}
			\noLine
			\UnaryInfC{$\sigma R n :: \Box A$}
		\end{prooftree}
				\ \\ \
		\begin{prooftree}
			\AxiomC{$\sigma :: \neg\Box_\rho A$}
			\LeftLabel{$(\pi_S, \Lambda)$}
			\RightLabel{, where $n\in\N$ is such that $\sigma S_\rho n\notin\Lambda$;}
			\UnaryInfC{$\sigma S_\rho n :: \neg A$}
			\noLine
			\UnaryInfC{$\sigma S_\rho n :: \Box A$}
		\end{prooftree}
				\ \\ \
		\begin{prooftree}
			\AxiomC{$\sigma :: \neg(A\vartriangleright B)$}
			\LeftLabel{$(\pi_\vartriangleright, \Lambda)$}
			\RightLabel{, where $n\in\N$ is such that $\sigma R n\notin \Lambda$.}
			\UnaryInfC{$\sigma R n :: A$}
			\noLine
			\UnaryInfC{$\sigma R n :: \Box_\sigma\neg B$}
			\noLine
			\UnaryInfC{$\sigma R n :: \Box\neg A$}
		\end{prooftree}
				\ \\ \
		We call the labelled formula above the line in the rules above the \emph{antecedent} and the labelled formula(s) under the line \emph{succedent(s)}.	
	\end{definition}
	
	Some clarifying remarks on the tableau rules seem in order. First we note that we use the symbol ``$\mid$" in the rules ($\rightarrow$), and ($\nu_\vartriangleright, \Lambda$) to denote branching in proof-trees. Next, we observe that various non-branching rules have multiple succedents such as the rules ($\neg \rightarrow$), $(\pi_\Box, \Lambda)$ and $(\pi_\vartriangleright, \Lambda)$.	These succedents are to be understood as different nodes one placed under the other. Lemma \ref{theorem:basicPropertiesIL}.\ref{theorem:basicPropertiesIL:Loeb} is reflected in the rules $(\pi_\Box, \Lambda)$ and $(\pi_\vartriangleright, \Lambda)$ and Lemma \ref{theorem:basicPropertiesIL}.\ref{theorem:basicPropertiesIL:LoebILversion} is reflected in the $(\pi_S, \Lambda)$ rule.
	
	Another non-local feature of the tableaux proof system will be that we will allow to apply rules to any node $\sigma :: A$ in a branch, not necessarily only to bottom-nodes. Upon application of the rule, the succedent(s) with possible branching can be appended to the bottom of \emph{any} branch passing through $\sigma:: A$. If $\B$ is a branch in a tree whose nodes are labelled formulas, by $\lab(\B)$ we denote the collection of labels that occur in $\B$.
	
	\begin{definition}[$\IL{\sf X}$-Tableaux, open and closed]\label{tableau}
		Given a Horn logic $\IL{\sf X}$ and a finite set $\Gamma$ of formulas, an \emph{$\IL{\sf X}$-tableau for $\Gamma$} is a binary irreflexive directed downward growing tree with nodes carrying labelled formulas defined inductively as follows:
		\begin{itemize}
			\item A single node tree $T$ with $0 :: A$ as the sole node for some formula $A\in \Gamma$ is an $\IL{\sf X}$-tableau for $\Gamma$.
			
			\item If $T$ is an $\IL{\sf X}$-tableau for $\Gamma$, then a tree $T'$ obtained by extending (appending below) any branches of $T$ with $0 :: A$ for some formula $A\in \Gamma$ is an $\IL{\sf X}$-tableau for $\Gamma$.
			
			\item Let $T$ be an $\IL{\sf X}$-tableau for $\Gamma$, $\B$ be a branch of $T$, and let $(\rho)$ be a rule w.r.t.~$\lab(\B)$.
			If some labelled formula $\sigma :: A$ that occurs in $\B$ is the antecedent of an instance of $(\rho)$, then the tree $T'$ obtained by extending $\B$ with the appropriate succedents of $(\rho )$ in any particular ordering (with possible branching) is an $\IL{\sf X}$-tableau for $\Gamma$.
		\end{itemize}
		
		A branch $\B$ of an $\IL{\sf X}$-tableau $T$ for $\Gamma$ is called \emph{closed} if there is $\sigma$ and $A$ such that $\sigma :: A\in\B$ and $\sigma :: \neg A\in\B$. Otherwise the branch is \emph{open}. An $\IL{\sf X}$-tableau $T$ for $\Gamma$ is closed if all of its branches are closed. Otherwise $T$ is open.
	\end{definition}
	
	Given a Horn logic $\IL{\sf X}$, we are now ready to assign to a finite set of formulas $\Gamma$ what we call a \emph{systematic $\IL{\sf X}$-tableau for $\Gamma$} which will contain all the information as to the satisfiability of $\Gamma$. The systematic tableau method given below follows closely the procedure given in \cite{Gore:1999}. 
	
	\begin{definition}[Systematic $\IL{\sf X}$-tableau]
		For a Horn logic $\IL{\sf X}$, a systematic $\IL{\sf X}$-tableau for a finite set $\Gamma$ of formulas is constructed in stages. Throughout the stages, the nodes in the tree $T_i$ will be marked with exactly one of \emph{awake, asleep} or \emph{finished}. The marked version of $T_i$ will be denoted by $\mu (T_i)$.
		
		\noindent\textbf{Stage 0:} Form \emph{the initial tableau} $T_0$ with $0 :: A$ for all $A\in \Gamma$ in some order on top of each other and mark them all awake.
		
		\noindent\textbf{Stage n+1:} Look for an awake $\sigma :: A$ in $\mu(T_n)$ closest to the root of the tableau; if there are several with the same distance, choose the leftmost one. If $A = p$ or $A = \neg p$ for some propositional variable $p$, then $T_{n+1}$ and $\mu(T_{n+1})$ are as $T_{n}$ and $\mu(T_{n})$ respectively except that we mark the node $\sigma :: A$ as finished and we end Stage n+1. 
		
		Otherwise we obtain $T_{n+1}$ and $\mu(T_{n+1})$ as follows:
		
		\begin{itemize}
			\item If $A = \neg\neg B$ for some $ B$, for every open branch $\B$ that passes through $\sigma :: A$, extend $\B$ with $\sigma ::  B$ marking it awake and marking $\sigma :: A$ as finished. Here and below `extending $\B$' means `appending new nodes to the bottom of $\B$'.
			
			\item If $A = ( B\rightarrow C)$ for some $ B$ and $ C$, for every open branch $\B$ that passes through $\sigma :: A$, split the end of $\B$ and extend the left fork with $\sigma :: \neg B$ and the right fork with $\sigma ::  C$. Both new nodes will be marked awake and $\sigma :: A$ will be marked as finished.
			
			\item If $A = \neg( B\rightarrow C)$ for some $ B$ and $ C$, for every open branch $\B$ that passes through $\sigma :: A$, extend $\B$ with $\sigma ::  B$ and $\sigma :: \neg C$ in whatever order. Both new nodes will be marked awake and $\sigma :: A$ will be marked as finished.
			
			\item If $A = \Box B$ for some $ B$, for every open branch $\B$ that passes through $\sigma :: A$ and for all $\tau\in\lab(\B)$, if $\sigma\R\tau$, then extend $\B$ with $\tau ::  B$. These new nodes will be marked awake and $\sigma :: A$ will be marked as asleep. 
			
			\item If $A = \neg\Box B$ for some $ B$, for every open branch $\B$ that passes through $\sigma :: A$, extend $\B$ with $\sigma R n :: \neg B$ and $\sigma R n :: \Box B$, where $n\in\N$ is the least number such that $\sigma R n\notin\lab(\B)$. Mark both $\sigma R n :: \neg B$ and $\sigma R n :: \Box B$ awake and $\sigma :: A$ finished. Moreover, mark as awake every $\tau :: \Box B\in\B$ and $\tau ::  B\vartriangleright C\in\B$ whenever $\tau\R\sigma R n$ and mark awake every $\tau :: \Box_\rho B\in\B$ whenever $\tau\S_\rho\sigma R n$. 
			
			\item If $A = \Box_\rho B$ for some $\rho$ and $ B$, for every branch $\B$ that passes through $\sigma :: A$ and for all $\tau\in\lab(\B)$, if $\sigma\S_\rho \tau$, then extend $\B$ with $\tau ::  B$. Mark $\tau ::  B$ awake and $\sigma :: A$ asleep. 
			
			\item If $A = \neg\Box_\rho B$ for some $\rho$ and $ B$, for every open branch $\B$ that passes through $\sigma :: A$, extend $\B$ with $\sigma S_\rho n :: \neg B$ and $\sigma S_\rho n :: \Box B$, where $n\in\N$ is the least number such that $\sigma S_\rho n\notin\lab(\B)$. Mark both $\sigma S_\rho n :: \neg B$ and $\sigma S_\rho n :: \Box B$ awake and mark $\sigma :: A$ finished. Moreover, mark awake every $\tau :: \Box B\in\B$ and $\tau ::  B\vartriangleright C\in\B$ such that $\tau\R\sigma S_\rho n$ and every $\tau :: \Box_\rho B\in\B$ such that $\tau\S_\rho\sigma S_\rho n$.
			
			\item If $A =  B\vartriangleright C$ for some $ B$ and $ C$, for every open branch $\B$ that passes through $\sigma :: A$ and every $\tau\in\lab(\B)$, if $\sigma\R\tau$, split the end of $\B$ and extend the left fork with $\tau :: \neg B$ and the right fork with $\tau :: \neg\Box_\sigma\neg C$. Both new nodes are marked awake and $\sigma :: A$ will be marked asleep.
			
			\item If $A = \neg( B\vartriangleright C)$ for some $ B$ and $ C$, for every open branch $\B$ that passes through $x$, pick the smallest $n{\in} \mathbb N$ such that $\sigma R n\notin \lab(\B)$ and extend $\B$ with $\sigma Rn ::  B$, with $\sigma Rn :: \Box_\sigma \neg C$ and with $\sigma Rn :: \Box\neg B$ in whatever order you like. All new nodes are marked awake and $\sigma :: A$ finished. Moreover, mark awake every $\tau :: \Box B\in\B$ and $\tau ::  B\vartriangleright C\in\B$ such that $\tau\R\sigma Rn$ and every $\tau :: \Box_\rho B\in\B$ such that $\tau\S_\rho\sigma Rn$. 
		\end{itemize}
		
		By this procedure we construct a chain $\langle T_i\colon i\in\omega\rangle$ of $\IL{\sf X}$-tableaux for $\Gamma$. We call $\bigcup_{i\in\omega} T_i$ \emph{a systematic $\IL{\sf X}$-tableau for $X$}.
	\end{definition}
	
	\begin{remark}
		A systematic $\IL{\sf X}$-tableau $T$ for a finite set $\Gamma$ of formulas is not in general\footnote{An example is when $\Gamma = \{  \Diamond p, p\rhd q, q\rhd p\}$: The $\Diamond p$ yields a world where $p$ holds, and since we do not reuse worlds the circular $p\rhd q, q\rhd p$ keep on creating new worlds.} an $\IL{\sf X}$-tableau in the sense of Definition \ref{tableau}. However, if $T$ is finite it is an $\IL{\sf X}$-tableau. In particular, if $T$ closes, then $T$ is an $\IL{\sf X}$-tableau. Moreover, if there is a closed $\IL$-tableau $T'$ for $\Gamma$, then there is a closed systematic tableau for $\Gamma$.
	\end{remark}
	
	\begin{lemma}[Fairness]
		If $\sigma :: A$ is awake at stage n+1, the systematic $\IL{\sf X}$-tableau procedure visits $\sigma :: A$ at some later stage.
		\begin{proof}
			Straightforward and similar to Lemma 6.4.4 in \cite{Gore:1999}.
		\end{proof}
	\end{lemma}
	
\section{Soundness}\label{section:soundness}
	As usual, for a Horn logic $\IL{\sf X}$, we call a formula $A$ $\IL{\sf X}$-tableau provable whenever the systematic $\IL{\sf X}$-tableaux for $\{  \neg A \}$ closes. In this section we shall show that this notion of provability is sound with respect to $\IL{\sf X}$ frames.
	
	\begin{definition}
		A set $\X$ of labelled formulas is $\IL{\sf X}$-satisfiable if there exists an $\IL{\sf X}$-model $M = \langle W,R,S,V\rangle$ and an \emph{interpretation} $I\colon\lab(\X)\rightarrow W$ such that
		\begin{itemize}
			\item[(i)] If $\sigma,\tau\in\lab(\X)$ and $\sigma\R\tau$, then $I(\sigma) R I(\tau)$;
			
			\item[(ii)] If $\rho,\sigma,\tau\in\lab(\X)$ and $\sigma\S_\rho\tau$, then $I(\sigma) S_{I(\rho)} I(\tau)$;
			
			\item[(iii)] $M,I(\sigma)\Vdash A$ for all $\sigma :: A\in\X$;
			
		\end{itemize}
		
		An $\IL{\sf X}$-tableau $T$ is $\IL{\sf X}$-satisfiable, if there is a branch $\B$ of $T$ such that $\B$ is $\IL{\sf X}$-satisfiable
	\end{definition}

The next lemma tells us that satisfiable $\IL{\sf X}$-tableaux are closed under applying the rules to them for Horn logics $\IL{\sf X}$.
 	
	\begin{lemma}
		Let $\IL{\sf X}$ be a Horn logic and let $T$ be a satisfiable $\IL{\sf X}$-tableau. Then for any rule, the tableau $T'$ obtained by the application of the rule is also $\IL{\sf X}$-satisfiable.
		\begin{proof}
			Suppose $\B$ is an $\IL{\sf X}$-satisfiable branch of $T$. We show that if we apply some rule to a some labelled formula in $\B$, we obtain a branch that is $\IL{\sf X}$-satisfiable. The cases for the {\bf propositional rules} are trivial.
			
			Suppose $\sigma :: \Box A\in\B$ and consider the branch $\C$ obtained by an {\bf application of the $(\nu_\Box)$-rule} with $\tau :: A$ added to the branch for some $\tau\in\lab(\B)$ such that $\sigma\R\tau$. By assumption there is an $\IL{\sf X}$-model $M = \langle W,R,S,V\rangle$ and an interpretation $I\colon\lab(\B)\rightarrow W$ such that $M,I(\sigma)\Vdash\Box A$. Now $\tau\in\lab(\B)$ and $\sigma\R\tau$. Hence $I(\sigma) R I(\tau)$ and so $M,I(\tau)\Vdash A$.
			
			Suppose $\sigma :: \Box_\rho A\in\B$ and consider the branch $\C$ obtained by an {\bf application of the $(\nu_S)$ rule} with $\tau :: A$ added to the branch for some $\tau\in\lab(\B)$ such that $\sigma\S_\rho \tau$. By assumption there is an $\IL{\sf X}$-model $M = \langle W,R,S,V\rangle$ and an interpretation $I\colon\lab(\B)\rightarrow W$ such that $M,I(\sigma)\Vdash\Box_{I(\rho)}A$. Now since $\tau,\rho\in\lab(\B)$ and $\sigma\S_\rho\tau$, we also have that $I(\sigma)S_{I(\rho)} I(\tau)$. Hence $M,I(\tau)\Vdash A$.
			
			Suppose $\sigma :: A\vartriangleright B\in\B$ and consider the two branches obtained by an {\bf application of the $(\nu_\vartriangleright)$ rule} with $\tau :: \neg A$ in the left branch and $\tau :: \neg\Box_\sigma\neg B$ in the right branch for some $\tau\in\lab(\B)$ such that $\sigma\R\tau$. Now by assumption there is an $\IL{\sf X}$-model $M = \langle W,R,S,V\rangle$ and an interpretation $I\colon\lab(\B)\rightarrow W$ such that $M,I(\sigma)\Vdash A\vartriangleright B$. Now since $\tau\in\lab(\B)$ and $\sigma\R\tau$, we have that $I(\sigma) R I(\tau)$. If $M,I(\tau)\Vdash \neg A$, then the left branch is satisfiable with $I$. If on the other hand $M,I(\tau)\Vdash A$, then there exists some $x\in W$ such that $I(\tau) S_{I(\sigma)} x$ and $M,x\Vdash B$. Hence $M,I(\tau)\Vdash\neg\Box_{I(\sigma)}\neg B$.
			
			Suppose $\sigma :: \neg\Box A\in\B$ and consider the branch $\C$ obtained by {\bf an application of the $(\pi_\Box)$ rule} with $\sigma Rn :: \neg A$ and $\sigma Rn :: \Box A$ added to the branch for some $n\in\N$ such that $\sigma R n\notin\lab(\B)$. By assumption there is an $\IL{\sf X}$-model $M = \langle W,R,S,V\rangle$ and an interpretation $I\colon\lab(\B)\rightarrow W$ such that $M,I(\sigma)\Vdash \neg\Box A$. Hence there is some $x\in W$ such that $I(\sigma) R x$ and $M,x\Vdash\neg A$ and $M,x\Vdash\Box A$. Now, since $\sigma Rn\notin\lab(\B)$, we can extend $I$ to $I'$ by putting $I'(\sigma Rn) = x$. Now define $\R_{I'}$ and $\S_{I'}$ on $\lab(\B)\cup\{\sigma Rn\}$ by
			\begin{align*}
			\langle \tau,\rho\rangle \in\R_{I'}&\Leftrightarrow I'(\tau)R I'(\rho),\\\\
			\langle \upsilon,\tau,\rho\rangle \in\S_{I'}&\Leftrightarrow I'(\tau) S_{I'(\upsilon)} I'(\rho).
			\end{align*}
			Now $\R_{I'}$ and $\S_{I'}$ satisfy conditions (\ref{definition:ILrelationsOnLabels:FirstItem})-(\ref{definition:ILrelationsOnLabels:LastItem}) in Definition \ref{definition:ILrelationsOnLabels}. Hence
			$\R\subseteq\R_{I'}$ and $\S\subseteq\S_{I'}$, and so $I'$ is an interpretation from $\C$ to $M$.
			
			Suppose $\sigma :: \neg\Box_\rho A\in\B$ and consider the branch $\C$ obtained by {\bf an application of the $(\pi_S)$ rule} with $\sigma S_\rho n :: \neg A$ and $\sigma S_\rho n :: \Box A$ added to the branch $\B$ for some $n\in\N$. By assumption there is an $\IL{\sf X}$-model $M = \langle W,R,S,V\rangle$ and an interpretation $I\colon\lab(\B)\rightarrow W$ such that $M,I(\sigma)\Vdash \neg\Box_{I(\rho)}A$. Now there is some $x\in W$ such that $I(\sigma) S_{I(\rho)} x$ and $M,x\Vdash\neg A$. Now if $M,x\Vdash\Box A$, we extend $I$ to $I'$ by putting $I'(\sigma S_\rho n) = x$. On the other hand if $M,x\Vdash\neg\Box A$, then there is $y\in W$ such that $x R y$ and $M,y\Vdash\neg A$ and $M,y\Vdash\Box A$. Now $I(\rho) R xRy$ and so $x S_{I(\rho)} y$. Hence $I(\sigma)S_{I(\rho)}y$. Now extend $I$ to $I'$ by putting $I'(\sigma S_{\rho} n) = y$. Now $I'$ is again an interpretation from $\C$ to $M$.
			
			Suppose finally $\sigma :: \neg(A\vartriangleright B)\in\B$ and consider the branch obtained by {\bf an application of the $(\pi_\vartriangleright)$ rule} with $\sigma R n :: A$, $\sigma Rn :: \Box_\sigma \neg  B$ and $\sigma Rn :: \Box\neg A$ added to the branch for some $n\in\N$. By assumption there is an $\IL{\sf X}$-model $M = \langle W,R,S,V\rangle$ and an interpretation $I\colon\lab(\B)\rightarrow W$ such that $M,I(\sigma)\Vdash \neg(A\vartriangleright B)$. Hence there is $x\in W$ such that $I(\sigma) R x$, $M,x\Vdash A$ and $M,x\Vdash \Box_{I(\sigma)}\neg B$. If $M,x\Vdash\Box\neg A$, we may extend $I$ to $I'$ by putting 
			$I'(\sigma Rn) = x$. On the other hand if $M,x\Vdash\neg\Box\neg A$, then there is $y\in W$ such that $M,y\Vdash A$ and $M,y\Vdash\Box\neg A$. But now since $I(\sigma) R x R y$, we have that $x S_{I(\sigma)} y$ and so $M,y\Vdash\Box_{I(\sigma)}\neg B$. So now we may extend $I$ to $I'$ by putting 
			$I'(\sigma R n) = y$. Again, $I'$ is an interpretation from $\C$ to $M$.
		\end{proof}
	\end{lemma}
	
	\begin{theorem}
Let $\IL{\sf X}$ be a Horn logic. If a systematic $\IL{\sf X}$-tableau for a set of formulas $\Gamma$ closes, then $\Gamma$ is $\IL{\sf X}$-unsatisfiable.
		\begin{proof}
			Suppose a systematic tableau $T$ for $\Gamma$ is closed, but that there is an $\IL{\sf X}$-model $M = \langle W,R,S,V\rangle$ and $x\in W$ such that $M,x\Vdash \Gamma$. Now consider the initial tableau $T_0$ for $\Gamma$ with $0 :: A$ for all $A\in \Gamma$. 
			
			Now, by assumption, $I = \{\langle 0,x\rangle \}$ is an interpretation from $\lab(T_0)$ to $W$. By the above lemma, every tableau obtained from the initial tableau is $\IL{\sf X}$-satisfiable. In particular, the closed tableau $T$ obtained by the systematic procedure is $\IL{\sf X}$-satisfiable. A contradiction since, for any branch $\B$ of $T$ there is $\sigma$ and $A$ such that $\sigma :: A\in\B$ and $\sigma :: \neg A\in\B$.
		\end{proof}
	\end{theorem}
For the sake of being explicit let us formulate the soundness of our tableaux as an immediate corollary.	
	\begin{corollary}
		Let $\IL{\sf X}$ be a Horn logic. If a systematic $\IL{\sf X}$-tableau for $\{\neg A\}$ closes, then $A$ is $\IL{\sf X}$-valid.
	\end{corollary}
	
\section{Completeness}\label{section:completeness}
	In this section we shall show that our proof system is also complete w.r.t.~$\IL{\sf X}$- frames.	
	\begin{definition}
		A set $\X$ of labelled formulas is a \emph{$\IL{\sf X}$-Hintikka set} if the following hold:
		\begin{itemize}
			\item[(i)] There is no $\sigma$ and $A$ such that $\sigma :: A\in\X$ and $\sigma :: \neg A\in\X$;\\
			
			\item[(ii)] If $\sigma :: \neg\neg A\in\X$, then $\sigma :: A\in\X$;\\
			
			\item[(iii)] If $\sigma :: A\rightarrow B\in\X$, then $\sigma :: \neg A\in\X$ or $\sigma ::  B\in\X$;\\
			
			\item[(iv)] If $\sigma :: \neg(A\rightarrow B)\in\X$, then $\sigma :: A\in\X$ and $\sigma :: \neg B\in\X$;\\
			
			\item[(v)] If $\sigma :: A\vartriangleright B\in\X$, then $\tau :: \neg A\in\X$ or $\tau :: \neg\Box_\sigma\neg B\in\X$ for all $\tau\in\lab(\X)$ such that $\sigma\R\tau$;\\
			
			\item[(vi)] If $\sigma :: \neg(A\vartriangleright B)\in\X$, then there is $\tau\in\lab(\X)$ such that $\tau :: A\in\X$, $\tau :: \Box_\sigma\neg B\in\X$ and $\sigma R\tau$;\\
		
			\item[(vii)] If $\sigma :: \Box A\in\X$, then $\tau :: A\in\X$ for each $\tau\in\lab(\X)$ such that $\sigma\R\tau$;\\
			
			\item[(viii)] If $\sigma :: \neg\Box A\in \X$, then there is $\tau\in\lab(\X)$ such that $\tau :: \neg A\in \X$ and $\sigma\R\tau$;\\
			\item[(ix)] If $\sigma :: \Box_\rho A\in\X$, then $\tau :: A\in\X$ for each $\tau\in\lab(\X)$ such that $\sigma\S_\rho\tau$;\\
			\item[(x)] If $\sigma :: \neg\Box_\rho A\in \X$, then there is $\tau\in\lab(\X)$ such that $\tau :: \neg A\in \X$ and $\sigma\S_\rho\tau$.\\
		\end{itemize}
	\end{definition}
Hintikka sets contain all the needed information to extract a model from them. This is clearly manifested in the proof of the following lemma.	
	\begin{lemma}[Truth Lemma]
Let $\IL{\sf X}$ be a Horn logic. If $\X$ is an $\IL{\sf X}$-Hintikka set and $\R^{\lab(\X)}_{\IL{\sf X}}$ is Noetherian, then $\X$ is $\IL{\sf X}$-satisfiable.
		
		\begin{proof}
			As mentioned before, we will omit various sub and superscripts. Thus, if $\R$ is Noetherian, then $\F = \langle \lab(\X),\R,\S\rangle$ is clearly an $\IL{\sf X}$-frame. Define a valuation $V$ on $\lab(\X)$ by putting
			$$V(p) = \{\sigma\in\lab(\X)\colon \sigma :: p\in\X \}\text{ for all propositional variables }p.$$
			Now we can prove by an easy induction on the complexity of formulas that for all $\sigma$ and $A$:
			\begin{align*}
			\sigma :: A\in\X &\Rightarrow \langle\F,V\rangle,\sigma\Vdash A\text{ and}\\
			\sigma :: \neg A\in\X &\Rightarrow \langle\F,V\rangle,\sigma\nVdash A.
			\end{align*}
			Hence $\langle \F,V\rangle$ satisfies $\X$ with the identity interpretation.
		\end{proof}
	\end{lemma}
	
	\begin{lemma}
		If $\B$ is an open branch in a systematic $\IL{\sf X}$-tableau for a finite set $\Gamma$, then $\B$ is a Hintikka set and $\R$ is Noetherian.
		\begin{proof}
			That $\B$ is a Hintikka set follows easily from the fairness of the systematic $\IL{\sf X}$-tableau procedure. 
			
			Notice that if $\sigma :: \Box A\in\B$ for some $\sigma$ and $A$, then $A$ is either a subformula of a formula from $\Gamma$ or the negation of a subformula of a formula from $\Gamma$. Let $\tilde{\Gamma} = \subjj[\Gamma]\cup\{\neg A\colon A\in\subjj[\Gamma]\}$. 
			
			Now suppose towards a contradiction that there is an ascending $\R$-chain $\langle \sigma_i\colon i\in\omega\rangle$ in $\lab(\B)$. Without loss of generality we may assume that $\sigma_0 = 0$ and $\sigma_i\neq 0$ for all $i > 0$.
			
			Now we show that for any $i\in\omega$ there is $A_i$ such that $\sigma_{i+1} :: \Box A_i\in\B$, but $\sigma_j :: \Box A_i\notin\B$ for all $j\leq i$. 
			
			If $\sigma_{i+1} = \tau R n$ for some $\tau\in\lab(\B)$ and $n\in\N$, then $\sigma_{i+1}$ is introduced either with a $(\pi_\Box)$-rule applied to some $\tau :: \neg\Box A\in\B$  or by a $(\pi_\vartriangleright)$-rule applied to some $\tau :: \neg(A\vartriangleright B)\in\B$. In the first case $\sigma_{i+1} :: \Box A\in\B$, but $\sigma_j :: \Box A\notin\B$ for all $j\leq i$. In the second case $\sigma_{i + 1} :: \Box\neg A\in\B$, but $\sigma_j :: \Box\neg A\notin\B$ for all $j\leq i$. 
			
			If $\sigma_{i+1} = \tau S_{\rho} n$ for some $\tau,\rho\in\lab(\B)$ and $n\in\N$, then $\sigma_{i+1}$ is introduced with a $(\pi_S)$-rule applied to some $\tau :: \neg\Box_\rho A$. Now $\sigma_{i+1} :: \Box A\in\B$, but $\sigma_j :: \Box A\notin\B$ for all $j\leq i$.
			
			Now for large enough $m$,
			$$|\{A\colon A\in\tilde{\Gamma}\text{ and }\sigma_m :: A\in\B\}| > |\tilde{\Gamma}|.$$ 
		\end{proof}
	\end{lemma}

We have the following corollaries from this lemma.
	
	\begin{corollary}
			Let $\IL{\sf X}$ be a Horn logic. If a systematic $\IL{\sf X}$-tableau for a finite $\Gamma$ has an open branch, then $\Gamma$ is $\IL{\sf X}$-satisfiable.
	\end{corollary}

In particular, we can formulate completeness of our proof systems. 
	
	\begin{corollary}
		Let $\IL{\sf X}$ be a Horn logic. If $A$ is $\IL{\sf X}$-valid, then any systematic $\IL{\sf X}$-tableau for $\{\neg A\}$ closes.
	\end{corollary}
	
\section{Scope of our results}\label{section:applications}
	
	Our results apply to all interpretability logics extending $\IL$ that are Horn. In particular the results apply to the most important systems $\IL\M$ and $\IL\P$. At first sight, the restriction of the logic being Horn might seem quite severe. However, most logics that occur in the literature turn out to be Horn. In particular also the logics based on 
	\begin{itemize}
		
		\item[$\sf R$:]
		$A \rhd B \ \to \ \neg (A \rhd \neg C) \rhd B\wedge \Box  C$;
	\end{itemize}
	and the two series of generalizations of this principle as presented in \cite{Joosten:2015:TwoSeries} are Horn logics. An important logic that falls out of the scope of this paper is $\IL{\sf W}$ since the corresponding frame condition is second order.

\section*{Acknowledgement}
We would like to thank Rajeev Gor\'e for his comments on a draft version of this paper. Further thanks go to Volodya Shavrukov and an anonymous referee for helping to improve the paper. The second author was supported by the Generalitat de Catalunya under grant number 2014 SGR 437 and from the Spanish Ministry of Science and Education under grant number MTM2014-59178-P.
	
	\bibliographystyle{plain}
	\bibliography{References}
\end{document}